\documentclass{article}

\usepackage{arxiv}

\usepackage[utf8]{inputenc} % allow utf-8 input
\usepackage[T1]{fontenc}    % use 8-bit T1 fonts
\usepackage{hyperref}       % hyperlinks
\usepackage{url}            % simple URL typesetting
\usepackage{booktabs}       % professional-quality tables
\usepackage{amsfonts}       % blackboard math symbols
\usepackage{nicefrac}       % compact symbols for 1/2, etc.
\usepackage{microtype}      % microtypography
\usepackage{lipsum}		% Can be removed after putting your text content
\usepackage{graphicx}

\usepackage{moreverb}

\usepackage{graphicx, subfigure}

\newcommand\BibTeX{{\rmfamily B\kern-.05em \textsc{i\kern-.025em b}\kern-.08em
T\kern-.1667em\lower.7ex\hbox{E}\kern-.125emX}}

\usepackage{graphicx}
\RequirePackage{fix-cm}
\usepackage{hyperref}
 \usepackage{epsfig} % for postscript graphics files
 \usepackage{amssymb}
 \usepackage{amsmath}
\usepackage{multirow}
\usepackage{epsfig} % for postscript graphics files
\input{epsf}
\usepackage{graphicx}
\usepackage{graphicx, subfigure}
\usepackage[utf8]{inputenc}
\usepackage[english]{babel}

\def\beq{\begin{equation}}
\def\eeq{\end{equation}}
\def\baq{\begin{eqnarray}}
\def\eaq{\end{eqnarray}}
\def\bal{\begin{align} }
\def\eal{\end{align} }
\def\bc{\begin{center}}
\def\ec{\end{center}}

\def\gr{\gamma_{r}}
\def\gr1{\gamma_{r1}}

\title{Development of New Methods for Detection \& Control of Credit Card Fraud Attacks  \\  }
\author{ Alexander Stotsky \\
Department of Computer Science and Engineering \\
     Chalmers University of Technology  \\
     Gothenburg  SE - 412 96, Sweden  \\
        \texttt{alexander.stotsky2@telia.com}  }
\hypersetup{
pdftitle={A template for the arxiv style},
pdfsubject={q-bio.NC, q-bio.QM},
pdfauthor={David S.~Hippocampus, Elias D.~Striatum},
pdfkeywords={First keyword, Second keyword, More},
}
\date{}
\begin{document}
\maketitle

\begin{abstract}
~Credit card fraud causes significant financial losses and frequently occurs as fraud attack, defined as short-term sequence of fraudulent transactions associated with high transaction rates and amounts, business areas historically tied to fraud, unusual transaction times and locations and different types of errors.
\\ Confidence interval method in the moving window with exponential forgetting is proposed in this report which allows to capture recent changes in the shopping behaviour of the cardholder, detect fraudulent amounts and mitigate the attack.
Fraud risk scoring method is used for estimation of the intensity of the fraudulent activity via monitoring of the transaction rates, merchant category codes, times and some other factors for detection of the start of the attack.
\\ The development and verification are  based on detailed analysis of the transaction patterns from the dataset, which represents an extensive collection of around $24.4$ million credit card transactions from IBM financial database.
Recommendations for further development of the detection techniques are also presented.
\end{abstract}


\maketitle

\section{Introduction}
\label{int}
\noindent
Automation of the payment processes in e-commerce reduces economic costs but
requires performance and robustness improvement of the real-time fraud detection systems
especially in the event of a crisis or war, \cite{riks}.
The  achievement of the trade-off between fraud detection accuracy  and  minimization of false alarms
(which is expensive and results in customer dissatisfaction) is the main challenge in real-time fraud detection systems.
\\
%Bayesian methods and neural networks are the most popular techniques in the credit card fraud detection.
Assuming that fraudulent transactions follow similar patterns the pattern recognition systems such as support vector machines, artificial neural networks, naive Bayesian networks, K - nearest neighbour methods, hidden markov models, fuzzy logic based systems and many other techniques can be applied for classification of transactions as fraudulent, see for example
\cite{pan} -\cite{haf} and references therein.
These methods do not take into account dynamics of fraud attacks (transient of fraudulent transactions) and
dynamics of  shopping behaviour of individual users  and based on average quantities of the large number of transactions,
which results in detection of almost every transaction as non-fraudulent due to overtraining and  wrong correlations, \cite{bani}.
\\ Time series analysis where each discrete step is associated with the transaction number, which allows identification of the dynamic spending behavior of individual cardholder and flagging transactions that differ from that behavior as outliers shows better outcomes.
\\
The aim of fraudsters, scammers and cybercriminals  is money and their behavior is associated with achievement of the highest benefits under time pressure.
Maximization of the returns results in sudden increases in fraudulent activity, which can be  defined as fraud attacks. Data analysis results in the following operational  definition of the fraud attack.
\\ {\it Fraud attack} is the sequence of fraudulent transactions (that can be conceptualized/framed as a dynamic system) which occur within the relatively short time interval  usually initiated by automated software tools (often coming from the same IP address) associated with high transaction rates, typical merchant category codes (MCC) (business areas traditionally
related to fraudulent activities), unusual transaction times and locations (often with geographical mismatches), technical glitches, insufficient balance, errors in CVV, PIN, zip codes and some other abnormal authorization errors.
\\
Outlier detection in time series with adjustable deviation (which takes into account the factors described above) from the
exponential moving average is the most efficient and implementable way for detection and control of fraud attack.
\\ Notice that the values of purchases in time series could be large for maximization of the outcome or low aiming for disguising fraud. Therefore the beginning of attack is usually associated with the sequence of low-value transactions (which makes the detection difficult), but the attack includes always the sudden changes in spending behaviour associated with high transaction amounts.  Detection of fraud attack can be more efficient than the detection of single fraudulent transaction due to availability of information about typical patterns/trends associated with elapsed time between transactions, merchant category codes, amounts, times, locations etc. Early detection of fraud attack that can be accomplished using these patterns is essential for minimizing financial losses and maintaining trust of the customer.
\\
Identification and description of the  purchasing profile of individual cardholder is absolutely necessary for improvement of the detection accuracy. Moreover, accurate detection is impossible without accounting for variability in purchasing behaviour of individual cardholder. The detection system should be adaptive taking into account recent changes in purchasing behaviour, which significantly reduces false alarms.
%Notice that introduction of such adaptive features in classical detection methods described above could be quite complicated.
Notice that inclusion of the adaptive features into classical detection techniques described above could pose considerable difficulties.
\\
Recent changes in the shopping profile of the cardholder  can be captured in a window of a certain size which is moving in time.  In addition to instantaneous forgetting of the shopping behaviour with moving window the exponential forgetting (where the weights decrease exponentially for more remote transactions) is introduced inside of the window.
In other words, combination of exponential forgetting and windowing techniques (where the forgetting factor and the window size are two adjustable parameters) opens new opportunities for performance improvement in time series analysis associated
with better trade-off achievements in the case of rapid changes in the shopping profile, \cite{jsp}. This results in significant improvement of the detection accuracy of fraud attacks.
\\
Data analysis shows that fraudulent activities usually occur as fraud attacks.
The shopping amounts are often not large in the beginning of the fraud attacks, which makes the detection difficult. However, detection of fraud attack can be more efficient than the detection of single fraudulent transaction due to availability of information about patterns of the attacks.
Therefore the fraud detection system should be based on combination of different features
which are able to correctly identify different parts of the fraud attack. The purchased amount should play the most important role in the detection of the attack. Taking into account that the amounts may be not large in the beginning of the attack  the inter-transaction time gaps, merchant category codes, transaction times and locations,  insufficient balance, technical glitches, errors in CVV, PIN, zip codes and other abnormal authorization errors should be used for identification of the start of the attack. The shopping amounts can then be controlled using confidence interval method minimizing losses.

\section{Confidence Interval for Detection \& Control}
\label{fdt}
\noindent
Outlier detection can be performed in moving window  and  associated with statistical hypothesis, where
weighted average is compared to suspected outlier within the two sample t-test framework, \cite{c}.
Namely, the hypothesis that the mean value is equal to the outlier is taken as a null hypothesis,
which is tested against the alternative hypothesis which indicates outlier in each step of sliding window.
Such formulation gives the basis for fraud  detection via the confidence interval method.
This confidence interval has two end points in each step of sliding window, which determine the fraud detection thresholds.
The endpoints are usually determined by the standard deviation, sample size (window size) and critical value taken from the Student distribution look-up table, under the assumption of normal distribution of the variable.
\\ Due to the significant differences and volatility in the shopping behavior of individual cardholders, assumptions about the distributions within the window cannot be made, and a problem-oriented decision regarding the length of the confidence interval as the detection threshold is needed.
In this particular problem the choice of the detection threshold is associated with the trade-off between fraud detection
accuracy and misdetection. The amounts associated with current shopping behaviour should be located
inside of the interval and all fraudulent amounts should exceed the threshold. It is clear that such trade-off is not
achievable without real-time adjustments of the window size, the forgetting factor and the length of the interval.
\\ The window size should be sufficiently small in order to capture fast changes in shopping behaviour
that limits application of different types of models. The accuracy of autoregressive moving average (ARMA) models, neural network models and some others is not sufficient for reliable detection of the outliers due to insufficient amount of data available inside  of short windows, \cite{mos}.
%Confidence interval method in moving window with prioritization of recent transactions is the only method for robust detection %of the outliers.
Confidence interval method based on the exponential moving average in short sliding window with prioritization of recent transactions is the only method for robust detection of the outliers in the case of rapid changes in shopping behavior.

\section{Detection of the Start of Attack: Fraud Risk Scoring Method}
\label{fdt1}
\noindent
Fraud risk scoring is well established technique widely used by financial institutions for detecting a range of fraudulent activities.
%Early detection of the start of the fraud attack is very important task for prevention of the financial losses.
The early detection of the start of the fraud attack is crucial for minimizing financial losses.
The start of the attack is usually characterized by relatively small amounts and the amount monitoring methods
cannot be used for detection, unfortunately.
\\
Therefore the detection of the start of the fraud attack should be based on estimation of the intensity of fraudulent activity.
{\it The intensity} refers to the degree of activity (associated with the frequency of fraudulent transactions) in the areas which are regularly linked to fraud.
\\
Fraud risk scoring is applied here for detection of the start of the fraud attack via estimation of the
intensity of the fraudulent activity. First, a number of factors associated with the strength of the attack
is identified for quantification.  Each factor in the list is quantified using risk score. The intensity is estimated
via total fraud risk score defined as the sum of individual scores of different factors.
The start of the attack is detected when the cumulative score goes beyond the predefined threshold.
%intensity of the
%activity in the areas which are frequently linked to fraud
%This method can also be used  for estimation of the likelihood of fraudulent transaction, where risk score is assigned
%to each event associated with the transaction. The decision is made based on the total fraud risk score defined as a sum of individual scores of different events. Risk scoring method is applied in this Section for detection of the start of the attack.
\begin{figure*}[!ht]
\begin{center}
  \includegraphics[width=16.61cm]{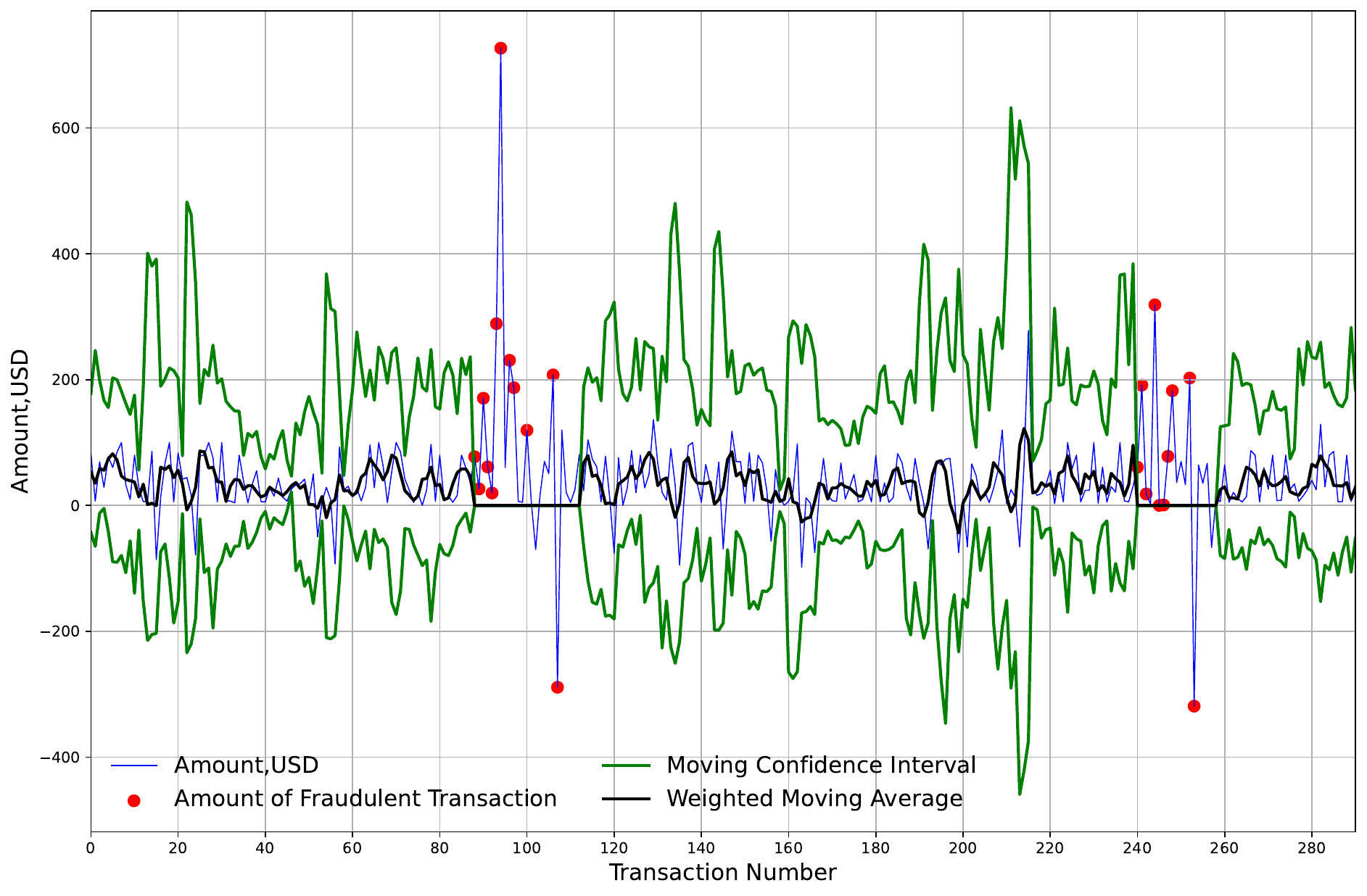}
\end{center}
\caption{\small{Detection and control of the fraud attack is presented in this Figure. The confidence is plotted with green lines.
The transaction amount is plotted with the blue line, weighted moving average is plotted with the black line
and fraudulent amounts are plotted  with red round signs.
Fraudulent transaction is detected if shopping amount exceeds the threshold. The first two transactions in the attack were flagged as fraudulent due small inter-transaction time gap, technical glitch and unusual time. The confidence interval was assigned to zero in each step of moving window and several online fraudulent transactions with high rates and amounts, which took place at unusual time  with technical glitches and CVV errors have been blocked and the fraud attack was prevented.
 \\ Similar procedure was applied to the second fraud attack shown in the Figure. The ability of the detection system to recover
 after the attacks is simply associated with resetting of the length of the confidence interval.
}}
\label{figci1}
\end{figure*}
\section{Control of Fraud Attack in the Case of Uncertainty}
\label{cont}
\noindent
When the start of the attack is detected the thresholds and the mean value are assigned to zero to control the  attack, see Figure~\ref{figci1}.  A number of online fraudulent transactions which took place at unusual times with high rates and amounts, well established MCC pattern, technical glitches and CVV errors has been blocked and the fraud attack was prevented. Notice that the average time of online attack is known and determines recovery after the attack, which is implemented as resetting of the length of confidence interval, see Figure~\ref{figci1}.
\\
Fraudsters continuously enhance their strategies and invent new forms of attacks  in emerging business areas such as quasi cash (cryptocurrency) and others. Existing detection algorithms may not be able to identify new patterns and attack mitigation strategies can be applied in the case of uncertainty, where the total score is not high, but the attack takes place.
Additional confirmation of the transaction in the form of single-use authentication code sent via sms can be requested in such cases. Data enrichment methods can also be applied in the case of uncertainty and additional relevant information such as details of the transaction, IP address, and other supporting information can be requested.
Finally, the transactions with certain merchant category codes which are likely associated with fraudulent activities
can be blocked or limited to small amounts, like $50$ USD, which results in significant savings.
%Notice that some other actions are also possible in the case of uncertainty (where total score is not high) in detection of the start of the attack.
\\ In the case of misdetection the transactions with merchant category codes which are rarely associated with fraudulent activities can be processed that increases customer satisfaction.
\\ Finally, the algorithm developers should monitor and proactively model new types of attacks which will likely appear in the future, outpacing fraudsters and protecting interests of the customers. To this end the simulation tools for modeling and detection of the attacks should be developed.
%Anticipating fraudsters' next moves
%Finally, the algorithm developers should monitor and proactively model new types of attacks which likely will appear in the future, predicting several steps ahead of the fraudsters and protecting interests of the customers. To this end the simulation tools for modeling and detection of the attacks should be developed, \cite{alt}.
\\ \\
Detailed description of detection algorithms is presented  in \cite{stata1} :
\\ \\
Stotsky A., Detection and Control of Credit Card Fraud Attacks in Sliding Window with Exponential Forgetting,
 International Journal of Computer Applications (0975 – 8887),
Volume 186 – No.74, March 2025, pp. 9-15.
 \\
     \url{https://www.ijcaonline.org/archives/volume186/number74/detection-and-control-of-credit-card-fraud-attacks-in-sliding-window-with-exponential-forgetting/}
\\ \\

\vspace{-8pt}
\section{Conclusion}
\label{con}
\noindent
New concept of fraud attack, which allows accurate detection and control the sequences of fraudulent transactions, reduces financial losses and increases the trust of customers was introduced in this report.
It was shown that the beginning of attack is usually associated with a number of small transactions and the attack includes always the sudden changes in spending behaviour associated with high transaction amounts. The detection of such attacks is a complicated  problem, which was divided in two parts.
\\ The first part is associated with monitoring of the transaction amounts using confidence interval method in the moving window with exponential forgetting which allows to capture recent changes in the shopping behaviour of the cardholder. This part is able to capture fraudulent transactions of relatively large amounts as deviations from current shopping
profile.
\\ The second part, which is assigned for detection of the beginning of the attack employs fraud risk scoring method for estimation of the intensity of fraudulent activity.
High transaction rates, business areas traditionally associated with fraudulent activities, unusual transaction times and locations, different types of errors and other factors are used  for detection of the start of the attack.
\\ The attack is controlled (after detection of the start) by adjusting the length of the confidence interval.
A number of control strategies can be used for mitigation of the attack: 1) all the transactions located outside of the interval can be blocked or limited to small amounts, 2) the transactions with merchant category codes which are likely associated with fraudulent activities can be blocked, 3) the transactions with merchant category codes which usually are not associated with fraudulent activities can be executed that increases customer satisfaction, and 4) many other actions.
\\ Further development of advanced real-time algorithms for detection of the fraud attacks based on data analysis and simulation tools is required since cybercriminals  continuously improve their strategies and create new methods. Developed detection algorithms should be robust against disturbances in  data communication systems and electricity grids for safe payments in crisis situations and war, \cite{riks}.

\section*{Disclosure Statement}
The dataset contains around $24.4$ million credit card transactions generated by  multi-agent virtual world simulation performed by IBM. The details of the generation method can be found in \cite{alt}. The simulations were performed in Python.

\end{document}